\documentclass[12pt]{article}
\usepackage{amssymb, amsmath, url, graphicx, setspace, geometry}
\usepackage{calrsfs}
\usepackage{wasysym}
\usepackage{xcolor}
\usepackage{tikz}

\def\3{\subset }
\def\4{\subseteq }
\def\<{\left<}
\def\>{\right>}

\def\bit{\begin{itemize}}
\def\eit{\end{itemize}}
\def\3{\subset }
\def\4{\subseteq }

\def\0{\leqno}

\def\barr{\begin{array}}
\def\earr{\end{array}}

\def\Z{{\rlap{$\kern2pt{\rm Z}$}{\rm Z}\,}}

\title{\bf Subgroups with a small sum of element orders}
\author{Mihai-Silviu Lazorec}
\date{August 2, 2022}

\begin{document}
\maketitle

\begin{abstract}
A finite group $G$ is said to be a $\mathcal{B}_{\psi}$-group if $\psi(H)<|G|$ for any proper subgroup $H$ of $G$, where $\psi(H)$ denotes the sum of element orders of $H$. In this paper, we characterize the $\mathcal{B}_{\psi}$-groups up to the class of finite nilpotent groups. We also study the above inequality in the case of finite $P$-groups and Schmidt groups, respectively. 
\end{abstract}



\section{Introduction}

Let $G$ be a finite group and let $n$ be a positive integer. We denote by $o(x), exp(G), \pi(G)$ and $C_n$ the order of an element $x\in G$, the exponent of $G$, the set of prime divisors of $|G|$ and the cyclic group of order $n$, respectively. By adding all element orders of $G$, we get the quantity
$$\psi(G)=\sum\limits_{x\in G}o(x)$$
known as the sum of element orders of $G$. During the past 12 years, this quantity was of some interest especially because it was proved that it can be used to characterize various classes of finite groups. More exactly, for a group $G$ of order $n$, it was shown that if $\frac{\psi(G)}{\psi(C_n)}>\alpha$, for some $\alpha\in (0, 1)$, then $G$ is cyclic/nilpotent/supersolvable/solvable (see \cite{3, 5, 6, 12, 17}). For each case, the best lower bound $\alpha$ was precisely determined as being $\alpha=\frac{\psi(H)}{\psi(C_{|H|})}$, where $H$ is the smallest  non-cyclic/non-nilpotent/non-supersolvable/non-solvable group. 

In \cite{1}, it was proved that the maximum value of the sum of element orders among all groups of order $n$ is $\psi(C_n)$. This means that if $G$ is a group of order $n$, then $\psi(G)\leq \psi(C_n)$ and the equality holds if and only if $G\cong C_n$. Obviously one may ask if something can be said about the minimum value of the sum of element orders. An answer is given up to nilpotent groups. More exactly, in \cite{2}, the authors proved the following result:\\

\textbf{Theorem 1.1.} \textit{Let $G$ be a group of order $n$. Then $\psi(G)\leq \psi(H)$ for every nilpotent group $H$ of order $n$ if and only if the Sylow subgroups of $G$ are of prime exponent.}\\

In \cite{11}, the authors classified the finite abelian groups determined by an inequality involving the sum of element orders. They said that a finite group $G$ is a \textit{$\mathcal{B}_{\psi}$-group} if the inequality 
$$\psi(H)< |G|$$
holds for all proper subgroups $H$ of $G$. Obviously, in order to inspect if $G$ is a $\mathcal{B}_{\psi}$-group, it suffices to check if the maximal subgroups of $G$ satisfy the above mentioned inequality. The classification of the finite abelian $\mathcal{B}_{\psi}$-groups is the following one:\\

\textbf{Theorem 1.2.} \textit{Let $G$ be a finite abelian group. Then $G$ is a $\mathcal{B}_{\psi}$-group if and only if $G\cong C_{p^2}$ or $G\cong C_p^n$, where $p$ is a prime and $n\geq 1$.}\\

Intuitively, both theorems conclude that if one wants to lower the sum of element orders of a finite group $G$, then a significant part of the group element orders should be equal to a prime/square-free number/square of a prime. The results of this paper also support this idea. Our main objective is to study the $\mathcal{B}_{\psi}$ property beyond finite abelian groups. A first step is to outline a connection between this property and a specific normal Hall $\pi$-subgroup of a finite group $G$. Then, we show that a finite nilpotent group $G$ is a $\mathcal{B}_{\psi}$-group if and only if $G\cong C_{p^2}$ or $G$ is a $p$-group with $exp(G)=p$. We enumerate the non-nilpotent solvable $\mathcal{B}_{\psi}$-groups up to order 2000. Such a group $G$ has the property $|\pi(G)|=2$. Inspired by this remark, we study the $\mathcal{B}_{\psi}$ property of two classes of finite groups whose orders are divisible by two primes: $P$-groups and Schmidt groups. 

We recall the definition of a finite $P$-group by following section 2.2 of \cite{15}. Let $p$ be a prime number and $n\geq 2$ be an integer. We say that a finite group $G$ belongs to the class $P(n, p)$ if $G$ is either an elementary abelian $p$-group of order $p^n$, or a semidirect product of an elementary abelian normal $p$-subgroup $H$ of order $p^{n-1}$ by a group of prime order $q\ne p$ which induces a non-trivial power automorphism on $H$. We call $G$ a $P$-group if $G$ belongs to $P(n, p)$, for some prime $p$ and some integer $n\geq 2$. If $p=2$, it is known that the elementary abelian 2-groups are the only finite $P$-groups, while if $p>2$, besides the elementary abelian $p$-group, one can construct a non-abelian $P$-group containing elements of order $q$, for every prime $q$ such that $q|(p-1)$.

We end this section by pointing out that some further research directions are outlined throughout the paper.

\section{Main results}

We begin this section by recalling some well-known results concerning the multiplicativity property of $\psi$ and the computation of the sum of element orders of a finite cyclic $p$-group. Also, we prove an inequality that is going to be used at some point.\\

\textbf{Lemma 2.1.} \textit{The following statements hold:
\begin{itemize}
\item[i)] $\psi(C_{p^n})=\frac{p^{2n+1}+1}{p+1}$, where $p$ is prime and $n\geq 1$ is an integer;
\item[ii)] Let $G_1$ and $G_2$ be finite groups. Then $\psi(G_1\times G_2)\leq \psi(G_1)\psi(G_2)$ and the equality holds if and only if $(|G_1|, |G_2|)=1$;
\item[iii)] Let $p, q$ be distinct primes and let $\alpha\geq 1, \beta\geq 2$ be some integers. Then $\psi(C_p^{\alpha})\psi(C_q^{\beta-1})> p^{\alpha}q^{\beta}.$
\end{itemize}}

\textbf{Proof.} For the first two items, the reader may check Lemma 2.1 of \cite{2} and Lemma 2.9 (1) of \cite{12}, respectively. For \textit{iii)}, we have
$$\psi(C_p^{\alpha})\psi(C_q^{\beta-1})> p^{\alpha}q^{\beta}\Longleftrightarrow (p^{\alpha+1}-p+1)(q^{\beta}-q+1)-p^{\alpha}q^{\beta}>0.$$
To prove the above inequality, it suffices to show that the function $f:[2,\infty)\longrightarrow \mathbb{R}$, given by $$f(x)=(px-p+1)(q^{\beta}-q+1)-xq^{\beta}, \ \forall \ x\in [2,\infty),$$ takes positive values. Obviously, $f$ is derivable and  
$$f'(x)=(p-1)q^{\beta}-pq+p\geq (p-1)q^2-pq+p>0,$$ so $f(x)\geq f(2), \ \forall \ x\in [2,\infty)$. Since $$f(2)=(p+1)(q^{\beta}-q+1)-2q^{\beta}\geq q^{\beta}-3q+3\geq q^2-3q+3>0,$$ it follows that $f$ takes only positive values. 
Hence, the proof is complete. 
\hfill\rule{1,5mm}{1,5mm}\\

Since a finite nilpotent group is isomorphic to the direct product of its Sylow subgroups, to investigate the $\mathcal{B}_{\psi}$ property of such a group, we need to solve the case of finite $p$-groups first.\\

\textbf{Proposition 2.2.} \textit{Let $G$ be a finite $p$-group. Then $G$ is a $\mathcal{B}_{\psi}$-group if and only if $G\cong C_{p^2}$ or $exp(G)=p$.}

\textbf{Proof.} Let $G$ be a $p$-group of order $p^n$, where $n\geq 1$. Assume that $G$ is a $\mathcal{B}_{\psi}$-group. If $n\in\lbrace 1,2\rbrace$, then $G\cong C_{p^2}$ or $exp(G)=p$, as desired. Let $n\geq 3$ and $exp(G)=p^m$, where $m\geq 1$. Clearly $m<n$ since otherwise $G$ would be a cyclic $\mathcal{B}_{\psi}$-group of order $|G|\geq p^3$. Then the conclusion of Theorem 1.2 would not hold, a contradiction. It follows that there is a maximal subgroup $H$ of $G$ containing an element of order $p^m$. We get
$$\psi(H)\geq 1+p^m+(p^{n-1}-2)p=|G|+p^m-2p+1.$$
If $m\geq 2$, we obtain $\psi(H)\geq |G|+(p-1)^2>|G|$, contradicting the $\mathcal{B}_{\psi}$ property of $G$. Hence, $m=1$ and $exp(G)=p$.

Conversely, if $G\cong C_{p^2}$, then $G$ is a $\mathcal{B}_{\psi}$-group by Theorem 1.2. If $exp(G)=p$, then for any maximal subgroup $H$ of $G$ we have
$$\psi(H)=1+(p^{n-1}-1)p=|G|-p+1<|G|,$$ so once again $G$ is a $\mathcal{B}_{\psi}$-group.   
\hfill\rule{1,5mm}{1,5mm}\\ 

In general, the direct product of two $\mathcal{B}_{\psi}$-groups is not a $\mathcal{B}_{\psi}$-group. For instance, according to Theorem 1.2 we know that $C_2$ and $C_3$ are $\mathcal{B}_{\psi}$-groups, while $C_6$ is not. Also, in general the $\mathcal{B}_{\psi}$ property of a finite group is not inherited by its non-trivial subgroups. The smallest example in this regard is $A_5$. Its maximal subgroups are isomorphic to $D_{10}, S_3$ or $A_4$. We have $\psi(D_{10})=\psi(A_4)=31$ and $\psi(S_3)=13$. All these sums of element orders are less than $|A_5|=60$, so $A_5$ is a $\mathcal{B}_{\psi}$-group. However, $\psi(C_5)=21>10=|D_{10}|$, so $D_{10}$ is not a $\mathcal{B}_{\psi}$-group. In what concerns the inheritance of the $\mathcal{B}_{\psi}$ property by non-trivial subgroups, we can use Proposition 2.2 to say more in the case of finite $p$-groups.\\

\textbf{Corollary 2.3.} \textit{Let $G$ be a finite $p$-group. Then all non-trivial subgroups of $G$ are $\mathcal{B}_{\psi}$-groups if and only if $G\cong C_{p^2}$ or $exp(G)=p$.}

\textbf{Proof.} Let $G$ be a finite $p$-group. If all non-trivial subgroups of $G$ are $\mathcal{B}_{\psi}$-groups, then $G$ is a $\mathcal{B}_{\psi}$-group. By Proposition 2.2, we deduce that $G\cong C_{p^2}$ or $exp(G)=p$. 

The converse easily follows if one uses Proposition 2.2 and the fact that $exp(G)=p$ implies that $exp(H)=p$ for any non-trivial subgroup $H$ of $G$.
\hfill\rule{1,5mm}{1,5mm}\\ 

Let $\pi$ be a set of prime numbers. We recall that a subgroup $H$ of a finite group $G$ is a Hall $\pi$-subgroup of $G$ if $\pi(H)\subseteq \pi$ and there is no prime divisor of $[G:H]$ which lies in $\pi$. In what follows, we show that if a finite group $G$ has a large normal Hall $\pi$-subgroup, then $G$ fails to be a $\mathcal{B}_{\psi}$-group.\\

\textbf{Proposition 2.4.} \textit{Let $G$ be a finite group and let $\pi(G)=\lbrace p_1, p_2, \ldots, p_k\rbrace$, where $k\geq 2$ and $p_1<p_2<\ldots<p_k.$ If $G$ has a normal Hall $\lbrace p_2, p_3, \ldots, p_k\rbrace$-subgroup, then $G$ is not a $\mathcal{B}_{\psi}$-group.}

\textbf{Proof.} Let $P_1$ be a Sylow $p_1$-subgroup of $G$ such that $|P_1|=p_1^{\alpha}$, where $\alpha\geq 1$. Let $H$ be a normal Hall $\lbrace p_2, p_3, \ldots, p_k\rbrace$-subgroup of $G$ and $K$ be a subgroup of $P_1$ of order $p_1^{\alpha-1}$. Since $H$ is normal, it follows that $HK$ is a maximal subgroup of $G$ of index $p_1$. Note that if $x\in H$ and $y\in HK\setminus H$, then $o(x)\geq p_2$ and $o(y)\geq p_1$. Then
$$\psi(HK)\geq 1+(|H|-1)p_2+|H|(|K|-1)p_1=|G|+(p_2-p_1)|H|+1-p_2>|G|$$
since $(p_2-p_1)|H|\geq p_2$. We conclude that $G$ is not a $\mathcal{B}_{\psi}$-group, as desired.
\hfill\rule{1,5mm}{1,5mm}\\

We have the necessary tools to generalize Theorem 1.2 up to the class of finite nilpotent groups.\\

\textbf{Theorem 2.5.} \textit{Let $G$ be a finite nilpotent group. Then $G$ is a $\mathcal{B}_{\psi}$-group if and only if $G\cong C_{p^2}$ or $exp(G)=p$, where $p$ is a prime.}

\textbf{Proof.} Let $G$ be a finite nilpotent $\mathcal{B}_{\psi}$-group. We take $\pi(G)=\lbrace p_1, p_2, \ldots, p_k\rbrace$, where $k\geq 1$, and denote by $P_i$ the Sylow $p_i$-subgroup of $G$, for all $i\in\lbrace 1,2,\ldots, k\rbrace$. Suppose that $k\geq 2$. Since $G$ is the direct product of its Sylow subgroups, we deduce that $P_2P_3\ldots P_k$ is a normal Hall $\lbrace p_2, p_3,\ldots, p_k\rbrace$-subgroup of $G$. According to Proposition 2.4, it follows that $G$ is not a $\mathcal{B}_{\psi}$-group, a contradiction. Then, $k=1$ and $G$ is a finite $p$-group satisfying the $\mathcal{B}_{\psi}$ property. By Proposition 2.2, we conclude that $G\cong C_{p^2}$ or $exp(G)=p$. 

As we previously saw, the converse also holds, so the proof is finished.
\hfill\rule{1,5mm}{1,5mm}\\

From now on, we provide some details concerning the non-nilpotent $\mathcal{B}_{\psi}$- groups. First of all, we note that Proposition 2.4 is also useful from a computational point of view, i.e. detecting examples of non-nilpotent $\mathcal{B}_{\psi}$-groups via GAP \cite{18}. For instance, there are 398,146,849 non-nilpotent groups of order $1536=2^9\cdot 3$. From these ones, 398,032,384 groups have a normal Sylow 3-subgroup, so they are not $\mathcal{B}_{\psi}$-groups by Proposition 2.4. The following table outlines the non-nilpotent $\mathcal{B}_{\psi}$-groups up to order 2000.\\

\begin{tabular}{ |p{3cm}|p{4cm}|p{2cm}|p{2cm}|  }
 \hline
 \multicolumn{4}{|c|}{Non-nilpotent $\mathcal{B}_{\psi}$-groups $G$ with $1\leq |G|\leq 2000$} \\
 \hline
 IdGroup & Structure Description & Solvable & $\pi(G)$\\
 \hline
 [12, 3]   &   $A_4$  & Yes & $\lbrace 2, 3\rbrace$ \\
 \hline
 [48, 50] & $C_2^4\rtimes C_3$  & Yes   & $\lbrace 2, 3\rbrace$ \\
 \hline 
 [56, 11] & $C_2^3\rtimes C_7$ & Yes & $\lbrace 2, 7\rbrace$\\
\hline
[60, 5] & $A_5$ & No & $\lbrace 2, 3, 5\rbrace$ \\
 \hline
 [80, 49] &  $C_2^4\rtimes C_5$  & Yes & $\lbrace 2, 5\rbrace$ \\
 \hline
 [120, 5] & $SL(2,5)$ & No &  $\lbrace 2, 3, 5\rbrace$\\
 \hline
 [160, 199] & $((C_2\times Q_8)\rtimes C_2)\rtimes C_5$  & Yes & $\lbrace 2, 5\rbrace$ \\
 \hline
 [168, 42] & $PSL(2,7)$ & No & $\lbrace 2, 3, 7\rbrace$\\
 \hline
 [192, 1541] & $C_2^6\rtimes C_3$ & Yes & $\lbrace 2, 3\rbrace$\\
 \hline
 [320, 1012] & $(C_2^2.C_2^4)\rtimes C_5$ & Yes & $\lbrace 2, 5\rbrace$\\
 \hline
 [336, 114] & $SL(2,7)$ & No & $\lbrace 2, 3, 7\rbrace$\\
 \hline
 [351, 12] & $C_3^3\rtimes C_{13}$ & Yes & $\lbrace 3, 13\rbrace$\\
 \hline
 [360, 118] & $A_6$ & No & $\lbrace 2, 3, 5\rbrace$\\
 \hline
 [405, 15] & $C_3^4\rtimes C_5$ & Yes & $\lbrace 3, 5\rbrace$\\
 \hline
 [448, 178] & $C_4^3\rtimes C_7$ & Yes & $\lbrace 2, 7\rbrace$\\
 \hline
 [448, 179] & $(C_2^3.C_2^3)\rtimes C_7$ & Yes & $\lbrace 2, 7\rbrace$\\
 \hline
 [448, 1393] & $C_2^6\rtimes C_7$ & Yes & $\lbrace 2, 7\rbrace$\\
 \hline
 [448, 1394] & $C_2^6\rtimes C_7$ & Yes & $\lbrace 2, 7\rbrace$\\
 \hline
 [504, 156] & $PSL(2,8)$ & No & $\lbrace 2, 3, 7\rbrace$\\
 \hline
 [576, 8661] & $C_2^6\rtimes C_9$ & Yes & $\lbrace 2, 3\rbrace$\\
 \hline
 [660, 13] & $PSL(2,11)$ & No & $\lbrace 2, 3, 5, 11\rbrace$\\
 \hline
 [720, 409] & $SL(2,9)$ & No & $\lbrace 2, 3, 5\rbrace$\\
 \hline
 [768, 1085321] & $C_2^8\rtimes C_3$ & Yes & $\lbrace 2, 3\rbrace$\\
 \hline
 [960, 11357] & $C_2^4\rtimes A_5$ & No & $\lbrace 2, 3, 5\rbrace$\\
 \hline
 [960, 11358] & $C_2^4\rtimes A_5$ & No & $\lbrace 2, 3, 5\rbrace$\\
 \hline
 [992, 194] & $C_2^5\rtimes C_{31}$ & Yes & $\lbrace 2, 31\rbrace$\\
 \hline
 [1092, 25] & $PSL(2,13)$ & No & $\lbrace 2, 3, 7, 13\rbrace$\\
 \hline
 [1215, 68] & $((C_3\times (C_3^2\rtimes C_3))\rtimes C_3)\rtimes C_5$ & Yes & $\lbrace 3, 5\rbrace$\\
 \hline
 [1280, 1116309] & $C_4^4\rtimes C_5$ & Yes & $\lbrace 2, 5\rbrace$\\
 \hline
 [1280, 1116310] & $((C_2\times((C_2\times((C_4\times C_2)\rtimes C_2))\rtimes C_2))\rtimes C_2)\rtimes C_5$ & Yes & $\lbrace 2, 5\rbrace$\\
 \hline
 [1280, 1116311] & $(C_2^4.C_2^4)\rtimes C_5$ & Yes & $\lbrace 2, 5\rbrace$\\
 \hline
 [1280, 1116312] & $((C_2^2\times Q_8)\rtimes Q_8)\rtimes C_5$ & Yes & $\lbrace 2, 5\rbrace$\\
 \hline
 [1280, 1116356] & $C_2^8\rtimes C_5$ & Yes & $\lbrace 2, 5\rbrace$\\
 \hline
\end{tabular} 
 
\begin{tabular}{ |p{3cm}|p{4cm}|p{2cm}|p{2cm}|  }
 \hline
 [1320, 13] & $SL(2,11)$ & No & $\lbrace 2, 3, 5, 11\rbrace$\\
 \hline
 [1344, 814] & $C_2^3.PSL(2,7)$ & No & $\lbrace 2, 3, 7\rbrace$\\
 \hline
 [1344, 11686] & $C_2^3\rtimes PSL(2,7)$ & No & $\lbrace 2, 3, 7\rbrace$\\
 \hline
 [1920, 240998] & $C_2^5\rtimes A_5$ & No & $\lbrace 2, 3, 5\rbrace$\\
 \hline
 [1920, 240999] & $C_2.(C_2^4\rtimes A_5)$ & No & $\lbrace 2, 3, 5\rbrace$\\
 \hline
 [1920, 241000] & $C_2^4\rtimes SL(2,5)$ & No & $\lbrace 2, 3, 5\rbrace$\\
 \hline
 [1920, 241001] & $C_2^5.A_5$ & No & $\lbrace 2, 3, 5\rbrace$\\
 \hline
 [1920, 241002] & $C_2^4\rtimes SL(2,5)$ & No & $\lbrace 2, 3, 5\rbrace$\\
 \hline
 \end{tabular}\\


We write some notes based on our enumeration. The list contains only three non-nilpotent $\mathcal{B}_{\psi}$-group of odd order. Some of the non-solvable examples are simple or quasisimple groups. In what concerns the non-nilpotent solvable examples, we observe a similarity: the orders of these groups have only two prime divisors. Hence, we pose the following open problem:\\

\textbf{Question 1.} \textit{Let $G$ be a finite non-nilpotent solvable $\mathcal{B}_{\psi}$-group. Is it true that $|\pi(G)|=2$?}\\

As a consequence of our last observation, one may be interested in studying the $\mathcal{B}_{\psi}$ property for specific classes of finite groups of order $p^{\alpha}q^{\beta}$, where $p, q$ are primes and $\alpha\geq 1, \beta\geq 1$. Two such remarkable classes are the finite non-abelian $P$-groups and the Schmidt groups (i.e. the finite minimal non-nilpotent groups). 

For fixed $p$ and $n$, a remarkable property of the finite $P$-groups is that they are lattice-isomorphic (see Theorem 2.2.3 of \cite{15}). This is an additional reason for which we choose to study the $\mathcal{B}_{\psi}$ property of the $P$-groups. More exactly, if we have two finite non-isomorphic groups $G_1$ and $G_2$ such that their subgroup lattices are isomorphic and $G_1$ is a $\mathcal{B}_{\psi}$-group, can we conclude that $G_2$ has the same property? The answer is negative since we prove that the only $P$-groups which satisfy the $\mathcal{B}_{\psi}$ property are the abelian ones.\\

\textbf{Corollary 2.6.} \textit{Let $G$ be a finite $P$-group. Then $G$ is a $\mathcal{B}_{\psi}$-group if and only if $G$ is elementary abelian.}

\textbf{Proof.} Let $G$ be a finite $P$-group. Suppose that $G$ is non-abelian and it satisfies the $\mathcal{B}_{\psi}$ property. By the description of the non-abelian $P$-groups, it follows that $G=H\langle x\rangle$, where $H$ is a normal Sylow $p$-subgroup of $G$ and $o(x)=q$.  Since $q|(p-1)$, it follows that $q<p$ so, according to Proposition 2.4, $G$ would not be a $\mathcal{B}_{\psi}$-group, a contradiction. Hence $G$ is an abelian $P$-group, so $G$ is elementary abelian, as desired. 

The proof is complete since the converse holds by Theorem 1.2. 
\hfill\rule{1,5mm}{1,5mm}\\

We move our attention to studying the $\mathcal{B}_{\psi}$ property of the Schmidt groups. These are the finite minimal non-nilpotent groups and  details on their properties and structure are covered in \cite{4, 14}. We recall some aspects as well. If $G$ is a Schmidt group, then $\pi(G)=\lbrace p, q \rbrace$, $G$ has a normal Sylow $p$-subgroup $P$ and a cyclic Sylow $q$-subgroup $Q$. We conclude that $G=P\rtimes Q$. Moreover, if $r$ is the multiplicative order of $p$ modulo $q$, it is known that if $P$ is abelian, then $P\cong C_p^r$, while if $P$ is non-abelian, then $Z(P)=P'=\Phi(P)$ is elementary abelian and $\frac{P}{\Phi(P)}\cong C_p^r.$ In short, $P$ is a special $p$-group (p. 183 of \cite{10}). Some constructions of special $p$-groups are covered in \cite{8}. Note that since $\Phi(P)$ and $\frac{P}{\Phi(P)}$ are elementary abelian, it follows that $exp(P)\in\lbrace p, p^2\rbrace$. Moreover, according to Theorem 5.2 of \cite{13}, if $p$ is odd, then $exp(P)=p$, while if $p=2$, then $exp(P)=4$. Finally, if $Q=\langle x\rangle$, then $x^q\in Z(G)$. Hence $\langle x^q\rangle$ is a normal subgroup of $G$ which trivially intersects $P$. Consequently, $G$ has a maximal subgroup $H\cong P\times\langle x^q\rangle$ of index $q$. All other maximal subgroups of $G$ are abelian and are of the form $\Phi(P)Q_i$ (see the proof of Theorem 1 of \cite{16}), where $i\in \lbrace 1, 2, \ldots, n_q\rbrace$ and $n_q$ is the number of Sylow $q$-subgroups of $G$. These would be the tools that are necessary to study the $\mathcal{B}_{\psi}$ property of the Schmidt groups.\\

\textbf{Theorem 2.7.} \textit{Let $G=P\rtimes Q$ be a Schmidt group of order $p^{\alpha}q^{\beta}$, where $p, q$ are distinct primes and $\alpha\geq 1, \beta\geq 1$. If $G$ is a $\mathcal{B}_{\psi}$-group, then $G=P\rtimes Q$, where $P$ is a special $p$-group of order $p^{\alpha}$, $Q\cong C_q$ and $p, q$ are primes such that $p<q$. Conversely,
\begin{itemize}
\item[i)] if $P$ is abelian, then $G$ is a $\mathcal{B}_{\psi}$-group;
\item[ii)] if $P$ is non-abelian, $p$ is odd, $r$ is the order of $p$ modulo $q$, and the inequality 
\begin{align}\label{r1}
(p^{\alpha-r+1}-p+1)(q^2-q+1)< p^{\alpha}q
\end{align}
holds, then $G$ is a $\mathcal{B}_{\psi}$-group;
\item[iii)] if $P$ is non-abelian, $p=2$, $r$ is the order of $p$ modulo $q$, and the inequalities (\ref{r1}) and
\begin{align}\label{r2}
\psi(P)<2^{\alpha}q
\end{align}
hold, then $G$ is a $\mathcal{B}_{\psi}$-group.
\end{itemize}
}

\textbf{Proof.} 
As stated above, $P$ is a normal Sylow $p$-subgroup isomorphic to a special $p$-group, $Q=\langle x\rangle\cong C_{q^{\beta}}$ and $H=P\times \langle x^q\rangle$ is a maximal subgroup of $G$. 

Assume that $G$ is a $\mathcal{B}_{\psi}$-group. Then, by Proposition 2.4, we have that $p<q$. Also, since $exp(P)\in\lbrace p, p^2\rbrace$, it follows that $\psi(P)\geq \psi(C_p^{\alpha})$. Using this fact and Lemma 2.1 \textit{i)} and \textit{ii)}, we have
$$\psi(H)=\psi(P)\psi(C_{q^{\beta-1}})\geq (p^{\alpha+1}-p+1)\cdot \frac{q^{2\beta-1}+1}{q+1}>p^{\alpha}\cdot \frac{q^{2\beta-1}+1}{q+1}.$$
If $\beta\geq 3$, then $\frac{q^{2\beta-1}+1}{q+1}>q^{\beta}$, so $\psi(H)>|G|$, contradicting the $\mathcal{B}_{\psi}$ property of $G$. Hence, $\beta\in\lbrace 1, 2\rbrace$. If $\beta=2$, according to Lemma 2.1 \textit{iii)}, we deduce that  
$$\psi(H)=\psi(P)\psi(C_q)\geq \psi(C_p^{\alpha})\psi(C_q)> p^{\alpha}q^2,$$ a contradiction. Consequently, $G=P\rtimes Q$ is a Schmidt group, where $P$ is a special $p$-group, $Q\cong C_q$ and $p<q$.

Conversely, let $G=P\rtimes Q$ be a Schmidt group of order $p^{\alpha}q$, where $\alpha\geq 1$, $P$ is a special $p$-group, $Q\cong C_q$ and $p, q$ are primes such that $p<q$. 

\textit{i)} If $P$ is abelian, then the maximal subgroups of $G$ are $P\cong C_p^{\alpha}$ and $Q_i\cong C_q$, where $i\in\lbrace 1, 2, \ldots, n_q\rbrace$. Also, $\alpha$ would be the multiplicative order of $p$ modulo $q$, so $q\leq p^{\alpha}-1$. We have $$\psi(P)=\psi(C_p^{\alpha})=p^{\alpha+1}-p+1<p^{\alpha+1}<p^{\alpha}q=|G|$$
and
$$\psi(Q_i)=\psi(C_q)=q^2-q+1<q^2\leq (p^{\alpha}-1)q<p^{\alpha}q=|G|, \ \forall \ i\in\lbrace 1, 2, \ldots, n_q\rbrace, $$
so $G$ is a $\mathcal{B}_{\psi}$-group. 

If $P$ is non-abelian, then the maximal subgroups of $G$ are $P$ and $\Phi(P)Q_i\cong C_p^{\alpha-r}\times C_q$, where $i\in\lbrace 1, 2, \ldots, n_q\rbrace$ and $r$ is the multiplicative order of $p$ modulo $q$. We distinguish the following two cases based on the parity of $p$.

\textit{ii)} If $p$ is odd, then $exp(P)=p$. Hence, $\psi(P)=p^{\alpha+1}-p+1<|G|$, as we previously saw. Assume that $G$ is not a $\mathcal{B}_{\psi}$-group. Then, there exists $i\in \lbrace 1, 2, \ldots, n_q\rbrace$ such that $$\psi(\Phi(P)Q_i)\geq |G|\Longleftrightarrow \psi(C_p^{\alpha-r})\psi(C_q)\geq |G|\Longleftrightarrow (p^{\alpha-r+1}-p+1)(q^2-q+1)\geq p^{\alpha}q,$$ a contradiction. Consequently, $G$ is a $\mathcal{B}_{\psi}$-group.

\textit{iii)} If $p=2$, then $exp(P)=4$. 
Assume that $G$ is not a $\mathcal{B}_{\psi}$-group. Then $\psi(P)\geq |G|$ or there exists $i\in \lbrace 1, 2, \ldots, n_q\rbrace$ such that $\psi(\Phi(P)Q_i)\geq |G|$. It follows that 
$$\psi(P)\geq 2^{\alpha}q \text{ or } (2^{\alpha-r+1}-1)(q^2-q+1)\geq 2^{\alpha}q,$$ 
a contradiction. Then $G$ is a $\mathcal{B}_{\psi}$-group and the proof is complete.
\hfill\rule{1,5mm}{1,5mm}\\

The result given by the case \textit{iii)} of the previous theorem can be improved if one explicitly computes $\psi(P)$. Since $P$ would be a special $2$-group with $exp(P)=4$ and $\Phi(P)\cong C_2^r$, we can easily get some bounds for $\psi(P)$ using the fact that the orders of the elements of $P\setminus \Phi(P)$ belong to $\lbrace 2,4\rbrace$. More exactly, we get
$$2^{\alpha+1}-1<\psi(P)\leq 1+2\cdot (2^{\alpha-r}-1)+4\cdot (2^{\alpha}-2^{\alpha-r}).$$ 

Also, concerning the cases \textit{ii)} and \textit{iii)} of the previous theorem, we note that there are examples of Schmidt groups for which at least one of the inequalities (\ref{r1}) and (\ref{r2}) hold or not. For instance, $SL(2, 3)\cong Q_8 \rtimes C_3$ is a Schmidt group having a non-abelian kernel $P\cong Q_8$ with $exp(P)=4$. This group is not a $\mathcal{B}_{\psi}$-group since (\ref{r2}) does not hold: we have $\psi(Q_8)=27>24=|SL(2,3)|$. Some examples of Schmidt $\mathcal{B}_{\psi}$-groups $G=P\rtimes Q$, such that $P$ is non-abelian and $exp(P)=4$, were outlined in our previous table: SmallGroup(160, 199), SmallGroup(320, 1012). Hence, both groups satisfy inequalities (\ref{r1}) and (\ref{r2}). Also, from the same table, we can extract SmallGroup(1215, 68) which is a Schmidt $\mathcal{B}_{\psi}$-group with non-abelian kernel of exponent $3$, so inequality (\ref{r1}) holds. Finally, it remains unknown if there are any Schmidt groups, with non-abelian kernel of odd exponent, which are not $\mathcal{B}_{\psi}$ groups. For instance, if one chooses the quadruple $(p, q, r, \alpha)=(3, 13, 3, 6)$, inequality (1) is not checked, but the only Schmidt group of order $9477$ is SmallGroup(9477, 505). However, its structure description is $C_{13}\rtimes C_{729}$, but we would need a non-abelian kernel of order 729 in order to follow the hypotheses of Theorem 2.7. Hence, we pose the following questions related to the last two paragraphs:\\

\textbf{Question 2.} \textit{Let $G$ be a non-abelian special $2$-group. What is the explicit form of $\psi(G)$?}\\

\textbf{Question 3.} \textit{Let $G=P\rtimes Q$ a Schmidt group of order $p^{\alpha}q$, where $\alpha\geq 1$ and $p, q$ are odd primes such that $p<q$. If $P$ is non-abelian, is it true that $G$ is a $\mathcal{B}_{\psi}$-group?}\\

If we check our previous table, we observe that $PSL(2, q)$ is a $\mathcal{B}_{\psi}$-group for $q\in \lbrace 3, 4, 5, 7, 8, 9, 11, 13 \rbrace$. For a prime power $q>3$, it is known that $PSL(2, q)$ is a non-abelian simple group. Hence, it would be interesting to check if the class consisting of these simple groups is contained in the class of $\mathcal{B}_{\psi}$-groups. For such a purpose, a main tool would be Dickson's classification of the maximal subgroups of $PSL(2, q)$ (see \cite{9}). One may ask:\\

\textbf{Question 4.} \textit{What can be said about the $\mathcal{B}_{\psi}$ property of the finite simple groups $PSL(2, q)$?}\\

Finally, other potential $\mathcal{B_{\psi}}$-groups are the $CP_1$ groups. The class of $CP_1$ groups consists of all finite groups whose element orders are primes. Hence, they are some strong candidates for our investigation. A complete classification of the $CP_1$ groups is highlighted in the main theorem of \cite{7}. Therefore, we suggest a final open problem:\\ 

\textbf{Question 5.} \textit{What can be said about the $\mathcal{B}_{\psi}$ property of $CP_1$ groups?}\\

A partial answer for this question is given in our paper since the nilpotent $CP_1$ groups are the finite $p$-groups of exponent $p$. These are $\mathcal{B}_{\psi}$-groups according to Proposition 2.2. The only non-solvable $CP_1$ group is $A_5$ which also is a $\mathcal{B}_{\psi}$-group. Hence, it remains to study the $\mathcal{B}_{\psi}$ property of the non-nilpotent solvable $CP_1$ groups. According to \cite{7}, such a group is a Frobenius group $G=P\rtimes Q$ of order $p^{\alpha}q$, where $p, q$ are primes and $\alpha\geq 1$. In addition, the kernel $P$ is of exponent $p$. 

\bigskip\noindent {\bf Acknowledgements.} The author is grateful to the reviewers for their remarks which improve the previous version of the paper. This work was supported by a grant of the "Alexandru Ioan Cuza" University of Iasi, within the Research Grants program, Grant UAIC, code GI-UAIC-2021-01.

\bigskip\noindent {\bf Conflict of interest.} The author have no conflicts of interest to declare. 

\bigskip\noindent {\bf Data availability.} Data sharing not applicable to this article as no datasets were generated or analysed during the current study.

\vspace*{3ex}
\small
\hfill
\begin{minipage}[t]{6cm}
Mihai-Silviu Lazorec \\
Faculty of  Mathematics \\
"Al.I. Cuza" University \\
Ia\c si, Romania \\
e-mail: {\tt silviu.lazorec@uaic.ro}
\end{minipage}

\end{document}